\documentclass{amsart}

\usepackage{amsfonts}
\usepackage{cases}
\usepackage{mathrsfs}
\usepackage{bbm}
\usepackage{amssymb}
\usepackage{txfonts}
\usepackage{amscd}
\usepackage{amsfonts,latexsym,amsmath,amsthm,amsxtra,mathdots}
\usepackage[bookmarks,colorlinks]{hyperref}
\usepackage[all,cmtip]{xy}
\RequirePackage{amsmath} \RequirePackage{amssymb}
\usepackage{color}
\usepackage{colordvi}
\usepackage{multicol}
\usepackage{tikz}
\usepackage{hyperref}
\usepackage{graphicx}
\usepackage{amsmath}
\usepackage{amsmath,amscd}
\usetikzlibrary{matrix,arrows}
\usepackage{textcomp}

     \newcommand{\BZ}{{\mathbb {Z}}}

\def\-{^{-1}}

\newcommand{\delete}[1]{}

    \theoremstyle{plain}

\newtheorem{thm}{Theorem}[section]
\newtheorem{lem}[thm]{Lemma}
\newtheorem{prop}[thm]{Proposition}

\newtheorem*{thmA}{Theorem A}

\newtheorem*{rem*}{Remark}

    \numberwithin{equation}{section}

\def\Proof{\noindent{\bf Proof}\quad}
\def\qed{\hfill$\square$\smallskip}

\begin{document}

\title{Farrell--Jones Conjecture for fundamental groups of graphs of virtually cyclic groups}

\author{Xiaolei Wu}
\address{Max Planck Institute for Mathematics, Bonn, Germany}
\email{hsiaolei.wu@mpim-bonn.mpg.de}

\subjclass[2010]{18F25,19A31,19B28}

\date{January, 2016}

\keywords{Farrell--Jones Conjecture; K-theory of group rings; L-theory of group rings; Bass--Serre theory.}

\begin{abstract}
In this note, we
prove the K- and L-theoretic Farrell--Jones Conjecture with coefficients in an additive category for  fundamental groups of  graphs of virtually cyclic groups.
\end{abstract}

\maketitle

\section*{Introduction}
In \cite{FW}, Farrell and the author proved the Farrell--Jones Conjecture for Baumslag-Solitar groups. Independently, Gandini--Meinert--R\"uping proved the Farrell--Jones Conjecture for any fundamental groups of graphs of abelian groups \cite{GMR}. This in particular includes all Baumslag--Solitar groups. One obvious question one can ask is whether the conjecture also holds for  fundamental groups of graphs of virtually abelian groups. One of the key technique difficulties that the previous methods do not work is that amalgamated products of virtually abelian groups in general do not have to be CAT(0); see \cite[III.$\Gamma$.6.13]{BH} for an example. However, amalgamated products of virtually cyclic groups are always CAT(0)  \cite[III.$\Gamma$.1.1 (c)]{BH}. In the present note, we are not able to deal with all the virtually abelian groups case. But we are able to show the following:

\begin{thmA}
The K- and L-theoretic Farrell--Jones Conjecture  with coefficients in an additive category and finite wreath products is true for fundamental groups of any graphs of virtually cyclic groups.
\end{thmA}
\begin{rem*} The problem whether  fundamental groups of any graphs of virtually cyclic groups satisfy the Farrell--Jones Conjecture was discussed before by Roushon in \cite{Ro}. He  obtained some partial results, for example, he showed  that fundamental groups of any trees of virtually cyclic groups satisfy the Farrell--Jones Conjecture \cite[Proposition 3.1.1]{Ro}.
\end{rem*}

Our theorem in particular implies any HNN extension of virtually cyclic groups satisfy the Farrell--Jones Conjecture. For more information about the Farrell--Jones Conjecture with finite wreath products and coefficients in an additive category we refer to \cite[Section 2.3]{W2}. In this note, we will abbreviate the K- and L-theoretic  Farrell--Jones Conjecture with finite wreath products and coefficients in an additive category by FJCw. The proof is similar to \cite{FW} and \cite{GMR} and depends heavily on them.  We try to map the group to some group that is known to satisfy FJCw, then use inheritance properties of FJCw to deduce the main result. The proof also uses crucially  a recent result of Bartles \cite[Remark 4.7, Example 4.10]{Ba} that FJCw is closed under taking
amalgamated products and HNN extensions along a finite subgroup.

\textbf{Acknowledgements.} The author wants to thank the Dahlem Research School at  Free University of Berlin and Max Planck Institute for Mathematics at Bonn for support. We also want to thank Daniel Kasprowski,  Henrik R\"uping, Sayed K. Roushon and F. Thomas Farrell for helpful comments.

\section{Inheritance Properties and Results on FJCw} \label{Inheritance-Properties}
We list some inheritance properties and results on FJCw that we may need.
\begin{prop}\label{qut}
(1) If a group $G$ satisfies FJCw, then any subgroup $H_1 \subset G$ and any finite index over group $H_2 \supset G$ satisfy FJCw.\\
(2) If $G_1$ and $G_2$ satisfy FJCw, then the direct product $G_1 \times G_2$ and the free product $G_1 \ast G_2$ satisfy FJCw.\\
(3) Let $\{G_i ~|~ i \in I\}$ be a directed system of groups (with not necessarily injective
structure maps). If each $G_i$ satisfies FJCw, then the colimit $colim_{i \in I} G_i$
satisfies FJCw.\\
(4) Let $\phi : G \rightarrow Q$ be a
group homomorphism. If $Q$, $Ker(\phi)$ and $\phi^{-1}(C)$ satisfy FJCw for every infinite cyclic subgroup $C < Q$ then G satisfies FJCw.\\
(5) CAT(0) groups satisfy FJCw.\\
(6) Virtually solvable groups satisfy FJCw.\\
(7) Fundamental groups of graphs of abelian groups satisfy FJCw.\\
(8) If $G_1$ and $G_2$ satisfy FJCw, then any amalgamated product of $G_1$ and $G_2$ along a finite subgroup and any HNN extension of $G_1$ along a finite subgroup satisfy FJCw.
\end{prop}

Proof of (1) - (4) can be found for example in \cite[Section 2.3]{W2} . (5) is the main result of \cite{BL} and \cite{W1} noticing that CAT(0) groups are closed under finite wreath products. (6) is proved in \cite{W2}. (7) is proved in \cite{GMR}. (8) is proved in \cite[Example 4.10]{Ba} for the case without finite wreath products and \cite[Remark 4.7]{Ba} explains how to extend it to the case with finite wreath products.

\section{Bass--Serre theory}
In this section, we review the basics of Bass--Serre theory. For more details, we refer to \cite[Chapter I]{DD} and \cite[Chapter 1]{Se}. We also prove a useful lemma.

Given a connected graph $\Gamma$ and an oriented edge $e\in E(\Gamma)$, denote its initial vertex by $\iota(e)\in V(\Gamma)$  and its terminal vertex by $\tau(e)\in V(\Gamma)$. Let $\overline{e}$ be the edge $e$ with opposite orientation, $\iota(\overline{e})=\tau(e)$ and $\tau(\overline{e})=\iota(e)$. A \emph{graph of groups structure} $\mathcal{G}$ on $\Gamma$ consists of families of groups $(G_v)_{v\in V(\Gamma)}$ and $(G_e)_{e\in E(\Gamma)}$ satisfying $G_{\overline{e}}=G_e$ for all $e\in E(\Gamma)$ and an injective group homomorphism $\alpha_e\colon G_e\hookrightarrow G_{\iota(e)}$ for each $e\in E(\Gamma)$. We call the pair $(\Gamma,\mathcal{G})$ a \emph{graph of groups}.

Given a maximal tree $T$ in $\Gamma$, let $\pi_1(\Gamma,\mathcal{G},T)$ be the group generated by the groups $G_v,v\in V(\Gamma)$ and the elements $e\in E(\Gamma)$ subject to the following relations
\begin{itemize}
\item[(i)] $\overline{e}=e^{-1}$ for all $e\in E(\Gamma)$;
\item[(ii)] $e\cdot\alpha_{\overline{e}}(s)\cdot\overline{e}=\alpha_e(s)$ for all $e\in E(\Gamma)$ and $s\in G_e$;
\item[(iii)] $e=1$ if $e\in E(T)$.
\end{itemize}
We call $\pi_1(\Gamma,\mathcal{G},T)$ the \emph{fundamental group} of $(\Gamma,\mathcal{G})$ relative to $T$. The isomorphism type of $\pi_1(\Gamma,\mathcal{G},T)$ does not depend on the choice of $T$, and we will just call it the \emph{fundamental group} of $(\Gamma,\mathcal{G})$ and denote it by $\pi_1(\Gamma,\mathcal{G})$. For each $v\in V(\Gamma)$ the canonical map $G_v\to\pi_1(\Gamma,\mathcal{G})$ is injective.  Given a graph of groups $(\Gamma,\mathcal{G})$, there is a corresponding Bass--Serre tree $X$ which $\pi_1(\Gamma,\mathcal{G})$ acts on with vertex stabilizer conjugating to some vertex group of $(\Gamma,\mathcal{G})$.

\begin{lem}\label{tree}
Let $(\Gamma,\mathcal{G})$ be a graph of groups and let $\phi: \pi_1(\Gamma,\mathcal{G}) \rightarrow \pi_1(\Gamma)$ be the quotient map induced by mapping all the vertex groups $G_v$ in $(\Gamma,\mathcal{G})$ to the trivial group. Let $X$ be the corresponding Bass--Serre tree of $(\Gamma,\mathcal{G})$. Then $ Ker(\phi) \backslash X$ is still a tree.
\end{lem}

\Proof
Note that $\pi_1(\Gamma, \mathcal{G})/Ker(\phi) \cong \pi_1(\Gamma)$ and thus we obtain an induced action of the free group  $\pi_1(\Gamma)$ on $Ker(\phi) \backslash X$. Since every vertex stabilizer of $\pi_1 (\Gamma,\mathcal{G})$ on $X$ is contained in $Ker(\phi)$, the action of $\pi_1(\Gamma)$ on $Ker(\Phi)\backslash X$ is free. We see that $Ker(\phi)\backslash X$ is the universal cover of $\Gamma$ and hence a tree.
\qed

\section{ Proof of Theorem A}
In this section, we prove our main theorem. We first prove FJCw for  fundamental groups of graphs of groups with vertex groups either $\BZ$ or $\BZ \rtimes \BZ/2$ and then use this to prove every fundamental group of a graph of virtually cyclic groups satisfies FJCw.

\begin{lem}\label{case:Z}
Let $(\Gamma, \mathcal{G})$ be a graph of groups with vertex groups either isomorphic to $\BZ$ or $\BZ\rtimes \BZ/2$, then $\pi_1(\Gamma, \mathcal{G})$ satisfies FJCw.
\end{lem}

\Proof By Proposition  \ref{qut} (3), FJCw is closed under directed colimit, we can assume the graph is finite. Proceeding by induction on the number of edges in $\mathcal{G}$, by Proposition \ref{qut}  (8), we only need to prove FJCw for graph of groups where all the edge groups are infinite. In fact, let $e \in \Gamma$ be an edge with edge group $G_e$ finite. If $\Gamma \setminus e $ is disconnected, let $\Gamma_1$ and $\Gamma_2$ be the two components. Restrict $(\Gamma, \mathcal{G})$  to the graphs $\Gamma_1$ and $\Gamma_2$, we get two graph of groups $(\Gamma_1, \mathcal{G}_1)$ and $(\Gamma_2, \mathcal{G}_2)$. Now $\pi_1(\Gamma, \mathcal{G})$ is a an amalgamated product of the group $\pi_1(\Gamma_1, \mathcal{G}_1)$ and $\pi_1(\Gamma_2, \mathcal{G}_2)$ along the edge group $G_e$. By Proposition \ref{qut}  (8), to prove FJCw for $\pi_1(\Gamma, \mathcal{G})$ it is enough to prove FJCw for $\pi_1(\Gamma_1, \mathcal{G}_1)$ and $\pi_1(\Gamma_2, \mathcal{G}_2)$. Similarly, if  $\Gamma \setminus e $ is connected, let $\Gamma'$ be the remaining graph and $(\Gamma', \mathcal{G}')$ be the corresponding graph of groups. Then $\pi_1(\Gamma, \mathcal{G})$ is an HNN extension of the group $\pi_1(\Gamma', \mathcal{G}')$ along the edge group $G_e$. Since $G_e$ is finite, by Proposition \ref{qut}  (8), to prove FJCw for $\pi_1(\Gamma, \mathcal{G})$ it is enough to prove FJCw for $\pi_1(\Gamma', \mathcal{G}')$. Inductively, we can get rid of all the finite edge groups in $(\Gamma, \mathcal{G})$ and reduce the proof of FJCw for $\pi_1(\Gamma, \mathcal{G})$ to the case where all the edge groups are infinite.

Now notice that any infinite subgroup of  $\BZ\rtimes \BZ/2$ is isomorphic to either $\BZ$ or  $\BZ\rtimes \BZ/2$, hence the edge groups of $(\Gamma, \mathcal{G})$ are either isomorphic to $\BZ$ or  $\BZ\rtimes \BZ/2$. We define another graph of groups $(\Gamma, \mathcal{G}')$ as follows. Note first that $\BZ \rtimes \BZ/2$ has an unique maximal infinite cyclic subgroup $\BZ\rtimes \{0\}$ and any infinite cyclic subgroup is contained in $\BZ\rtimes \{0\}$. Thus for each vertex group $G_v$, we can define $G'_v$ to be the quotient of $G_v$ by its maximal infinite cyclic subgroup, hence $G'_v$ equals $\BZ/2$ if $G_v$ is isomorphic to $\BZ\rtimes \BZ/2$, $\{0\}$ otherwise. We define the edge group $G'_e$ similarly. For each edge $e$, denote its initial vertex by $\iota(e)$ and its terminal vertex by $\tau(e)$. Corresponding to each edge $e$, there is an injective group homomorphism $\alpha_e: G_e \hookrightarrow G_{\iota(e)}$. Our quotienting out the maximal infinite cyclic subgroup process also induces a map $\alpha'_e: G'_e \hookrightarrow G'_{\iota(e)}$ and we claim that  it is well defined and still injecitive, thus we complete our definition of  $(\Gamma, \mathcal{G}')$. Here is the proof of this claim. If $G_e$ is isomorphic to $\BZ$, then after quotienting out the maximal infinite cyclic subgroup, $G_e$ becomes the trivial group. Hence $\alpha'_e$ is well defined and injective. If $G_e$ is isomorphic to $\BZ\rtimes \BZ/2$, then after quotienting out maximal infinite cyclic subgroups, both $G_e $ and $ G_{\iota(e)}$ becomes $\BZ/2$. Notice that $\BZ\rtimes \{0\} \subset G_e$ must be mapped to a subgroup of $\BZ\rtimes \{0\} \subset G_{\iota(e)}$, thus $\alpha'_e$ is well defined. Notice also that the order $2$ element $(0,1)\in \BZ\rtimes \BZ/2 \cong G_e$ must be mapped to an element of the form $(a,1) \in  \BZ\rtimes \BZ/2 \cong G_{\iota(e)}$, thus $\alpha'_e$ is still injective.

 Since its vertex groups are either trivial or $\BZ/2$,   $\pi_1(\Gamma, \mathcal{G}')$ act on its corresponding Bass--Serre tree properly and cocompactly. Thus $\pi_1(\Gamma, \mathcal{G}')$ is a CAT(0) group and satisfies FJCw. Moreover, there is an obvious map $q: \pi_1(\Gamma, \mathcal{G})\rightarrow \pi_1(\Gamma, \mathcal{G}')$ which is induced by quotienting out  the maximal infinite cyclic subgroup at each vertex group. By Proposition \ref{qut} (4), we are left to verify  $Ker(q)$ and $q^{-1}(C)$ satisfy FJCw for any infinite cyclic subgroup $C$ of $\pi_1(\Gamma, \mathcal{G}')$. By Bass--Serre theory there is a tree $X$ that $\pi_1(\Gamma, \mathcal{G})$ acts on with vertex stabilizers isomorphic to some conjugate of the vertex groups. Now restricted to each vertex stabilizer, $Ker(q)$ becomes $\BZ$. Hence $Ker(q)$ is the fundamental group of a graph of infinite cyclic groups. By Proposition \ref{qut} (7), $Ker(q)$ satisfies FJCw. Now for $q^{-1}(C)$, since each vertex group in $(\Gamma, \mathcal{G}')$ is trivial or $\BZ/2$ hence $C$ does not lie in any vertex group or any conjugates of it. Thus $q^{-1}(C)$ acts on $X$ with infinite cyclic stabilizers and satisfies FJCw by Proposition \ref{qut} (7). \qed

To proceed we need the following lemma of Farrell--Jones \cite[Lemma 2.5]{FJ}.

\begin{lem}\label{vcstru}
Let $G$ be a virtually cyclic group. Then $G$ contains a unique maximum normal finite subgroup $F$ with one of the following holds

\begin{itemize}

\item (1) the finite case, $G =F$;
\item (2) the orientable case, $G/F$ is the infinite cyclic group;
\item (3) the nonorientable case, $G/F$ is the infinite dihedral group.

\end{itemize}

\end{lem}

We proceed to prove our main theorem. Let $(\Gamma, \mathcal{G})$ be a graph of groups with vertex groups virtually cyclic. By the arguments at the first paragraph of the proof of Lemma \ref{case:Z}, we can assume the group $\Gamma$ is finite and all the edge groups are infinite.

At each vertex group $G_v$, let $F_v$ be its unique maximum normal subgroup determined by Lemma \ref{vcstru}. We define a new graph of groups $(\Gamma, \mathcal{G})$ as follows. For each vertex $v \in \Gamma$, let $G'_v := G_v / F_v$. Note that $G'_v$ is isomorphic to either $\BZ$ or $\BZ \rtimes \BZ/2$. For each edge $e$, let the corresponding injective group homomorphism be $\alpha_e\colon G_e\hookrightarrow G_{\iota(e)}$. Define $G'_e = G_e/\alpha_e^{-1}(F_{\iota(e)})$, denote the induced edge map by $\alpha'_e$. We need to show that $G'_e = G'_{\bar e}$ and $\alpha'_e$ is still injective. Note that $\alpha_e^{-1}(F_{\iota(e)})$ is a finite normal subgroup, we only need to show that $\alpha_e^{-1}(F_{\iota(e)})$ is  the unique finite maximal normal subgroup of $G_e$ in Lemma \ref{vcstru}. Let $x$ be an element of infinite order in $G_e$, then $\alpha_e(x) \not \in F_{\iota(e)}$ since $\alpha_e$ is injective. Hence $[x] \in G_e/\alpha_e^{-1}(F_{\iota(e)}) $ and its image $\alpha'_e([x]) \in G_{\iota(e)}  /F_{\iota(e)}$ also have infinite order. In particular the image $\alpha'_e(G_e/\alpha_e^{-1}(F_{\iota(e)})) \subset G_{\iota(e)}  /F_{\iota(e)}$ is an infinite group. Since $G_{\iota(e)}  /F_{\iota(e)}$ is isomorphic to $\BZ$ or $\BZ \rtimes \BZ/2$, $\alpha'_e(G_e/\alpha_e^{-1}(F_{\iota(e)})) $ is also isomorphic to $\BZ$ or $\BZ \rtimes \BZ/2$. Let $q_e$ be the quotient map from $G_e$ to $G_e/\alpha_e^{-1}(F_{\iota(e)})$, then the kernel of the composite $\alpha'_e \circ q_e: G_e \hookrightarrow G_{\iota(e)}/ F_{\iota(e)}$ is the unique finite normal subgroup by Lemma \ref{vcstru}. On the other hand the kernel of $\alpha'_e \circ q_e$ is precisely  $\alpha_e^{-1}(F_{\iota(e)})$. This  shows  $\alpha_e^{-1}(F_{\iota(e)})$ is the unique finite maximal normal subgroup and also the induced map $\alpha'_e\colon G'_e\hookrightarrow G'_{\iota(e)}$ is injective.

So far, we have defined a new graph of groups  $(\Gamma, \mathcal{G}')$ and there is an obvious map $q: \pi_1 (\Gamma, \mathcal{G}) \rightarrow  \pi_1(\Gamma, \mathcal{G}')$ induced by quotienting out the unique maximal normal finite subgroup at each vertex group. By Lemma \ref{case:Z}, $\pi_1(\Gamma, \mathcal{G}')$ satisfies FJCw. By Proposition \ref{qut} (4), we are left to verify  $Ker(q)$ and $q^{-1}(C)$ satisfy FJCw for every infinite cyclic subgroup $C$ of $\pi_1(\Gamma, \mathcal{G}')$. Again, by Bass--Serre theory there is a tree $X$ that $\pi_1(\Gamma, \mathcal{G})$ acts on with vertex stabilizers isomorphic to some conjugate of the vertex groups. As a subgroup of $\pi_1(\Gamma, \mathcal{G})$, $Ker(q)$ also acts on $X$. Now restricted to each vertex stabilizer, $Ker(q)$ is finite. Hence $Ker(q)$ is the fundamental group of a graph of finite groups. These groups are directed colimit of CAT(0) groups since they act on their corresponding Bass--Serre trees with finite stabilizers,  therefore $Ker(q)$ satisfies FJCw; see for example \cite[Lemma 3.3]{GMR} for more details. Now for an infinite cyclic group $C$, if it intersects any
conjugate of a vertex group in $(\Gamma, \mathcal{G}')$ trivially, then $q^{-1}(C)$ also acts on $X$ with finite stabilizers and  satisfies FJCw. On the other hand, if $C$ intersects some conjugate of a vertex group nontrivially, up to finite index, we can assume $C$ lies in some conjugate of  a vertex group by Proposition \ref{qut} (1). As in Lemma \ref{tree}, we can map $\pi_1(\Gamma, \mathcal{G}')$ to $\pi_1(\Gamma)$ defined by mapping every vertex group to the trivial group, denote the map by $\phi$. Now $q^{-1}(C)$ is a subgroup of $Ker(\phi \circ q)$. On the other hand, we can apply Lemma \ref{tree} to the map $\phi\circ q$ and we see that $Ker(\phi \circ q) \backslash X$ is a tree. Moreover, $Ker(\phi \circ q)$ acts on $X$ with virtually cyclic stabilizers. By \cite[III. $\Gamma$.1.1 (c)]{BH} any amalgamated products of two CAT(0) groups along a virtually cyclic group is again a CAT(0) group. Hence $Ker(\phi \circ q)$ is a directed colimit of CAT(0) groups and satisfies FJCw. \qed

\section{Further remarks}

In general, let $(\Gamma, \mathcal{G})$ be a graph of virtually abelian groups. By  Proposition  \ref{qut} (3), we can assume the graph is finite and each vertex group is finite generated. Now given a virtually abelian group $A$, it has a finite index free abelian subgroup $B$ of finite rank. we can embed $A$ into a bigger group $B  \wr F$ for some finite group $F$(see for example the proof in \cite[Proposition 2.17]{W2}). Note that $B \wr F = B^{|F|} \rtimes F$ further embeds in $\BZ \wr S_m = \BZ^m \rtimes S_m$ where $m= rank(B) |F|$ and $S_m$ is the symmetric group of $m$ elements. Since the graph has only finite many edges, we can further embed every vertex group into $\BZ^{M} \rtimes S_{M}$ for some $M$ sufficient large. Now using this, we can embed $\pi_1(\Gamma, \mathcal{G})$ into a graph of groups where every vertex group is $\BZ^{M} \rtimes S_{M}$. Since FJCw is closed under taking subgroups (Proposition \ref{qut} (1) ), we have proved the following:

\begin{prop}
If for any  $m$, FJCw holds for fundamental groups of graphs of groups with vertex groups isomorphic to $\BZ^{m} \rtimes S_{m}$, then FJCw holds for fundamental groups of graphs of virtually abelian groups.
\end{prop}

This seems to reduce the problem to  easier cases. However, we do not know how to map these groups to some nice groups that are known to satisfy FJCw.

\bibliographystyle{amsplain}

\end{document}